\documentclass[12pt]{article}

\usepackage{latexsym}
\usepackage{amsopn,amscd}
\usepackage{amsfonts}
\usepackage{amssymb}
\usepackage{amsthm}
\usepackage{amsmath}
\usepackage{a4wide}
\newtheorem{Theorem}{Theorem}[section]
\newtheorem{Proposition}[Theorem]{Proposition}

\newtheorem{Remark}[Theorem]{Remark}
\newtheorem{Example}[Theorem]{Example}

\def\SSS{S}
\def\AA{A}
\numberwithin{equation}{section}

\begin{document}

\title{\bf Braids and crossed modules}

\author{
J.~Huebschmann
\\[0.3cm]
 USTL, UFR de Math\'ematiques\\
CNRS-UMR 8524
\\
59655 Villeneuve d'Ascq C\'edex, France\\
Johannes.Huebschmann@math.univ-lille1.fr
 }
\maketitle
\begin{abstract}
{For $n \geq 4$, we describe Artin's braid group on $n$ strings as
a crossed module over itself. In particular, we interpret the braid
relations as crossed module structure relations.}
\end{abstract}

\noindent
\\[0.3cm]
{\bf 2010 Subject classification:}~ 
{Primary: 20F36 57M05 57M20 57M25; Secondary: 18D50 18G55 20C08 55P48}
\\[0.3cm]

\section{Introduction}

The  braid group $B_n$ on $n$ strings  was
introduced by  Emil Artin  in 1926,
cf. \cite{magnuone} and the references there. 
It has various    interpretations, specifically, 
as  the group of geometric braids in 
${\mathbb R}^3$, as the mapping class group of
an $n$-punctured   disk, etc.
The group $B_1$ is the trivial group and $B_2$ is infinite cyclic.
Henceforth we will take $n \geq 3$.
The group $B_n$ has $n-1$ generators
$\sigma_1$,\ldots , $\sigma_{n-1}$ subject to the relations
\begin{align}
\sigma_i\sigma_j&=\sigma_j\sigma_i,\ 
\vert i-j\vert \geq 2,
\label{rel01}
\\
\sigma_j\sigma_{j+1} \sigma_j&= \sigma_{j+1} \sigma_j\sigma_{j+1},
\ 1 \leq j \leq n-2.
\label{rel02}
\end{align}
The main result of the present paper, Theorem \ref{second} below, 
says that,
as a crossed module over itself,
(i) 
the Artin braid group $B_n$ has a single generator, which can be taken to be 
any of the Artin generators
$\sigma_1$,\ldots , $\sigma_{n-1}$, and that,
furthermore, (ii) the kernel of the surjection
from the free $B_n$-crossed module $C_n$ 
in any one of the $\sigma_j$'s
onto $B_n$ coincides with
the second homology group $\mathrm H_2(B_n)$ of $B_n$, well known to be cyclic 
of order 2 when $n \geq 4$ and trivial for $n=2$ and $n=3$.
Thus,
the crossed module property corresponds precisely to the relations
\eqref{rel02}, that is, these relations amount to {\em Peiffer identities\/}
(see below).
When
$C_n$ is taken to be the
free $B_n$-crossed module in $\sigma_1$, for $n \geq 4$, the element 
${}^{\sigma_3}\sigma_1\sigma_1^{-1}$ of $C_n$ has order 2 and generates the copy of $\mathrm H_2(B_n)\cong \mathbb Z/2$ in $C_n$;
the element ${}^{\sigma_3}\sigma_1\sigma_1^{-1}$ of $C_n$ being the generator
of the kernel of the projection from $C_n$ to $B_n$, it is then manifest how
this central element forces the relations \eqref{rel01}.

The braid groups are special cases of the more general class
of Artin groups \cite{briesait}, \cite{deligone}.
However, for a general Artin group presentation, there does not seem
to be an obvious way to concoct a crossed
module of the same kind as that for the braid group presentation,
and the results of the present paper 
seem to be special for braid groups.

More recently augmented racks \cite{fennrour} 
and quandles have been explored in the literature, and there is an intimate
relationship between augmented racks, quandles, crossed modules, 
braids, tangles, etc.
\cite{fennrour}.
See also Remark \ref{crack} below.
Is there a single theory, yet to be uncovered,
so that augmented racks, quandles, links, braids, tangles, 
crossed modules, Hecke algebras, etc.
are different incarnations thereof?
In Section \ref{geometrici} below we explain how operads are behind 
our approach and, perhaps, operads can be used to build such a single theory.
We hope our paper contributes to the general understanding of the situation.

It is a pleasure to acknowledge the stimulus of conversation
with L. Kaufman at a physics meeting in Varna (2007)
which made me rethink  old ideas of mine related with braids,
identities among relations,
and crossed modules.
I am indebted to
J. Stasheff for a number of comments on a draft of the paper
which helped improve
the exposition and to F. Cohen for discussion about the generating cycle
of $\mathrm H_2(B_n)$, see Remark \ref{fcoh} below; this discussion 
prompted me to work out the relationship between polyhedral geometry 
and identities among relations in Section \ref{geometrici} below.

\section{Historical remarks}

The defining relations \eqref{rel01} and \eqref{rel02}
for $B_n$ are intimately connected with  
Hecke algebra relations.   
V. Jones discovered
new invariants of knots \cite{vjoneone}
by analyzing certain finite-dimensional 
representations of the Artin braid group
into the Temperley-Lieb algebra \cite{templieb},
that is, a von Neumann algebra
defined by a kind of Hecke algebra relations.

In the late 1970's, 
working on the unsettled problem of J. H. C. Whitehead's
\cite{whitethr} whether any subcomplex of an aspherical
2-complex is itself aspherical
\cite{aspheric}, I had noticed that the relations \eqref{rel02}
strongly ressemble the {\em Peiffer identities\/}
(cf. \cite{peiffone}, \cite{reideone})
defining a free crossed module,
but it seems that this  
similarity has never been made explicit in the literature
nor has it been used to derive structural insight
into braid groups or into related groups such as knot 
or link groups.

That kind of  similarity gets even more striking when
one examines Whitehead's original proof \cite{whitethr}
of the following result of his: 
{\em Adding 2-cells to a pathwise connected space $K$ so that the space
$L$ results 
yields the $\pi_1(K)$-crossed
module $\pi_2(L,K)$ that is {\em free\/} on the 2-cells that have been attached
to form $L$.\/}
A precise form of this result will be spelled out below as Theorem 
\ref{whit}.
The crucial step of that proof proceeds as follows:
Let  $B^2$ denote the 2-disk, $S^1$ its boundary circle,
and denote by $o$ the various base points.
To prove that the obvious map from the free
$\pi_1(K)$-crossed module on the 2-cells in $\mathrm L \setminus K$ to
$\pi_2(L,K)$ is injective,
Whitehead considers a  null-homotopy 
$(B^2\times I, S^1\times I,o\times I) \to (L,K,o)$;
when this null-homotopy is made transverse, the pre-images of suitable
small disks in the interiors of the 
2-cells in $\mathrm L \setminus K$ form
a link in the solid cylinder $B^2\times I$
(more precisely: a partial link based at the 2-disk $B^2$ at the bottom
of the cylinder)
which truly arises from a kind of tangle or braid;
Whitehead then argues in terms of the exterior of that link
and in particular interprets the Wirtinger relations corresponding
to the crossings in terms of the identities defining a free crossed module,
that is, in terms of the Peiffer identities.
In the late 1970's I learnt from R. Peiffer (private communication)
that Reidemeister had also been aware of the appearance of this kind of
link resulting from a null homotopy.
A more recent account of Whitehead's argument can be found in
\cite{browwhit}.
See also the survey article \cite{browhueb}, where
the relationship between identities among relations and
links is explained in detail, in particular in terms of the
notion of {\em picture\/} which is a geometric representation
of an element of a crossed module; a geometric object equivalent
to a picture had been exploited already in \cite{peiffone}
under the name {\em Randwegaggregat\/}.
These links very much look like the tangles that have been used
to describe the Temperley-Lieb algebra, cf. \cite{kauflins} (p.~8) 
which, in turn, leads, via a trace, to a description of the Jones polynomial. The formal relationship between Peiffer identities and
this kind of link was also exploited in \cite{sieraone}.

In the early 1980's, some of my ideas related with identities among relations
and crossed modules
merged into the papers
\cite{cch}, \cite{collhueb}, and \cite{aspheric};
in particular, pictures are crucially used in these papers.

\section{Crossed modules}

We will denote the identity element of a group by $e$.
Let $C$ and $G$ be groups, $\partial \colon C \to G$ 
a homomorphism,
and suppose that the group $C$ is endowed with an action of $G$
which we write as
\[
G \times C \longrightarrow C,\ (x,c) \longmapsto {}^xc,\ x \in G,\ x \in C.
\]
These data constitute a {\em crossed module\/} 
provided 
\begin{itemize}
\item
relative to the conjugation action of $G$ on itself
$\partial$ is a homomorphism of $G$-groups, and 
\item
in $C$
the identities
\begin{equation}
aba^{-1} = {}^{\partial a}b
\label{peiffer}
\end{equation}
are satisfied.
\end{itemize}
We refer to the identities \eqref{peiffer}
or, equivalently,
\begin{equation}
aba^{-1}  ({}^{\partial a}b)^{-1}=e,
\label{peiffer3}
\end{equation}
as {\em Peiffer identities\/}.
Morphisms of crossed modules are defined in the obvious way.
Thus crossed modules form a category.

Let $\partial \colon C \to G$
be a crossed module.
For later reference, we recall the following facts:
\begin{itemize}
\item
The image  $\partial (C)\subseteq G$ is a normal subgroup; we will write
the quotient group as $Q$.
\item
The kernel $\mathrm{ker}(\partial)$ of $\partial$ is a central subgroup of $C$
and acquires a $Q$-module structure.
\item
Via the $G$-action, the abelianized group $C_{\mathrm{Ab}}$
acquires a $Q$-module structure.
\end{itemize}

Let $\varphi \colon Y \to C$ be a map.
The crossed module $\partial \colon C \to G$
is said to be {\em free\/} on $\varphi$ provided the following 
universal property is satisfied:

Given a crossed module $\delta \colon D \to H$,
a homomorphism $\alpha \colon G \to H$,
and a (set) map $b \colon Y \to D$  such that
$\alpha\circ \partial\circ \varphi = \delta \circ b$, 
there is a unique group homomorphism
$\beta \colon C \to D$ with $\beta\circ \varphi =b$ such that
the diagram
\begin{equation}
\begin{CD}
C @>{\partial}>> G
\\
@V{\beta}VV
@VV{\alpha}V
\\
D @>{\delta}>> H
\end{CD}
\end{equation}
is a morphism of crossed modules from $\partial$ to $\delta$.

When  $\partial \colon C \to G$
is free on $\varphi\colon Y \to C$, 
in view of the universal property of a free crossed module,
the abelianized group
$C_{\mathrm{Ab}}$ is then the free $Q$-module having (a set in bijection with)
$Y$ as a $Q$-basis.
Consequently, the map $\varphi$ is then injective and,
somewhat loosely, one  says that 
 $\partial \colon C \to G$
is {\em free\/} on $Y$.

Given the crossed module $\partial \colon C \to G$,
we will refer to $C$ as a $G$-{\em crossed module\/};
when  $\partial \colon C \to G$
is a free crossed module,
we will refer to $C$ as a {\em free\/} $G$-crossed module.

Let $L$ be a space obtained from the pathwise connected
space $K$ by the operation
of attaching 2-cells.
In \cite{whitethr}, J.H.C. Whitehead attempted an algebraic description
of the second homotopy group $\pi_2(L)$. In \cite{whitefou},
he reformulated his earlier results in terms of a precise algebraic
description of the second relative homotopy group $\pi_2(L,K)$.
In \cite{whitefiv}, he finally codified the result 
as follows:

\begin{Theorem}[Whitehead]
\label{whit}
Via the homotopy boundary map  from
$\pi_2(L,K)$ to $\pi_1(K)$, the group
$\pi_2(L,K)$ is the free crossed $\pi_1(K)$-module in the 2-cells
that have been attached to $K$ to form $L$.
\end{Theorem}

Recall that the {\em geometric realization\/} 
of a presentation $\langle X,R\rangle$ of a group is the 2-complex
$K(X,R)$ with a single vertex whose 1- and 2-cells are
in bijection with $X$ and $R$, respectively,
in such a way that, for $r \in R$, 
the attaching map of the associated 2-cell corresponds to the word
$\xi_r$ in the generators $X$ when the relator $r$ is written out as a word
in the generators $X$.

Let $\langle X,R\rangle $ 
be a presentation of a group, let $F$ be the free group on $X$,
and let $N$ be the normal closure in $F$ of the relators $R$.
An {\em identity among relations\/} \cite{peiffone}
(synonymously: identity sequence) 
is an $m$-tuple $(q_1,\ldots,q_m)$ ($m$ may vary)
of pairs of the kind $q_j= (\xi_j, r_j^{\epsilon_j})$
where $\xi_j \in F$, $r_j \in R$, $\epsilon_j = \pm 1$, such that
evaluation of $\prod \xi_j r_j^{\epsilon} \xi_j^{-1}$ in $F$ yields
the identity element $e$ of $F$. Among the observations 
made in \cite{peiffone}, we quote only the following ones:

\noindent
(i) Given two relators $r$ and $s$,
the quadruple
\begin{equation}
((e,r), (e,s), (e,r^{-1}), (\xi_r, s^{-1})),
\label{peiffer2}
\end{equation}
where $\xi_r \in F$ refers to the evaluation of $r$ in $F$,
is always an identity among relations, independently of a particular
presentation. In the context of group presentations,
these identities are referred to
in the literature as {\em Peiffer identities\/}.

\noindent
(ii) The group of identities among relations modulo that of identities
of the kind \eqref{peiffer2}
is isomorphic to the second homotopy group
of the geometric realization of the presentation $\langle X,R\rangle$.
Here the term \lq\lq second homotopy group\rq\rq\ is to be 
interpreted in the language of homotopy chains, which here captures the
second homology group of the universal covering space of the geometric
realization of the presentation under discussion.

\noindent
(iii) The combinatorial structure
of a 3-manifold can be characterized in terms of identities among relations;
in particular, the operation of attaching a 3-cell
admits an interpretation in terms of killing
a certain identity among relations. 

In \cite{reideone}, some of these results
are rediscussed and the idea of a free crossed module 
is introduced without a name.
A modern account of identities among relations and crossed modules
can be found in \cite{browhueb}.

\begin{Example} 
{\rm Given the group $G$, with respect to the adjoint action,
the identity homomorphism $G \to G$
is a crossed module in an obvious manner.
In the present paper, we will explore this crossed module for the special case
where $G$ is the braid group $B_n$ on $n$ strings.}
\end{Example}

\begin{Example} 
{\rm The injection of a normal subgroup into a group
is a crossed module in an obvious way.}
\end{Example}

\begin{Example} 
{\rm Given a group $G$, the obvious map $G \to \mathrm{Aut}(G)$ which sends
a member of $G$ to the inner automorphism determined by it is a crossed 
module. This observation has been explored in \cite{crossed}.}
\end{Example}

\begin{Example}\label{ex2} 
{\rm Given the central extension
$1
\to
\AA
\to
E
\to
G
\to
1
$
of (discrete) groups, where $\AA$ is an abelian group,
conjugation in $E$ induces a crossed module structure on $E \to G$.}
\end{Example}

\begin{Example} {\rm Let $\langle X,R\rangle$ be the Wirtinger presentation
of a knot in $\mathbb R^3$, let $K$ be the geometric realization
of $\langle X,R\rangle$, pick a generator $x \in X$, and let $r$ be the
relator $r=x$. The normal closure of $r$ is well known to be the
entire knot group.
Let $ L= K\cup e^2$ be the geometric realization
of the resulting  presentation
$\langle X, R\cup\{r\}\rangle$ of the trivial group.
By construction, $\pi_1(K)$ is the group of the knot
and, in view of Whitehead's
theorem (Theorem \ref{whit} above),  the homotopy boundary map
$\partial \colon \pi_2( L, K) \to  \pi_1(K)$
is a free crossed module, 
the $\pi_1(K)$-crossed module $\pi_2(L, K)$
being free on the single generator $[r]$.
Since $L$ is actually contractible,
$\partial \colon \pi_2( L, K) \to  \pi_1(K)$
is an isomorphism.
Thus the group of the knot is
a free crossed module over itself in a single generator.
Suitably interpreted, this kind of observation is also
valid for
a geometrically indecomposable link.}
\end{Example}

\begin{Example}\label{rack} {\rm Let $G$ be a group. A $G$-{\em rack\/}
or {\em augmented rack\/} {\rm \cite{fennrour}}
or $G$-{\em automorphic set\/} {\rm \cite{briesthr}}
is a $G$-set $Y$ together with a map $\partial \colon Y \to G$
such that 
\begin{equation}
\partial ({}^xy) = x\partial(y) x^{-1},  \ x \in G,\ y \in Y.
\label{rackpeiffer}
\end{equation}
The universal object associated with the $G$-rack $Y$ is a crossed module
$\partial\colon C_Y\to G$ 
{\em \cite{fennrour}}
where the notation $\partial$ is abused. In particular, when
$Y$ is the free $G$-set $G \times S$ for a set $S$,
the crossed module
$\partial\colon C_Y\to G$ is the free crossed module on $S$.
Notice that the identities {\rm \eqref{rackpeiffer}}
yield again a version of the Peiffer identities.}
\end{Example}

\section{The main technical result}

\begin{Theorem}\label{basic} 
Let $G$ be a group, let $\SSS\subseteq G$
be a subset of $G$ whose normal closure is all of $G$, 
let $\partial \colon C \to G$ be the free crossed module on $\SSS$,
and suppose that the induced map $C_{\mathrm{Ab}} \to G_{\mathrm{Ab}}$
is an isomorphism of abelian groups.
Then $\partial$ fits into a central extension of the kind
\begin{equation}
\begin{CD}
0 @>>> \mathrm H_2(G) @>>> C @>{ \partial}>> G @>>> 1.
\end{CD}
\label{central8}
\end{equation}
\end{Theorem}

\begin{proof}
Let $\langle X, R\rangle$ be a presentation of $G$,
let $F$ be the free group on $X$, and let $N$ be the normal closure
of $R$ in $F$. Abusing notation,
let $\SSS\subseteq F$ be a subset which under the projection
from $F$ to $G$ is mapped bijectively onto $\SSS \subseteq G$.
Then $\langle X, R\cup S\rangle$ is a presentation of the trivial group.

Let 
$K$ be the geometric 
realization of the presentation
$\langle X, R\rangle$ and $L$ 
that of the presentation $\langle X, R\cup S\rangle$
of the trivial group.
By construction, the fundamental group $\pi_1(K)$
of $K$ is $G$. 
In view of Whitehead's
theorem 
(Theorem \ref{whit} above),  the homotopy boundary map
$\partial \colon \pi_2(L,K) \to  G$
is a free crossed module, 
the $G$-crossed module $C=\pi_2(L,K)$
being free on the set $\SSS$. 

Since $L$ is simply connected,
the Hurewicz map $\pi_2(L)\to \mathrm H_2(L)$ is an isomorphism
and
the exact homotopy sequence of the pair $(L,K)$ 
takes the form
\begin{equation}
\begin{CD}
\pi_2( K) @>>> \mathrm H_2(L) @>>> \pi_2(L, K) @>{ \partial}>> G @>>> 1.
\label{central22}
\end{CD}
\end{equation}

It is a classical fact that the cokernel of the Hurewicz map
$\pi_2(K) \to \mathrm H_2(K)$ is the second homology group $\mathrm H_2(G)$
of $G$ \cite{hopfone}.
The group $\mathrm H_2(L,K)$ is free abelian on a basis in bijection with
the set $S$, and
inspection of the exact homology sequence
\begin{equation}
\begin{CD}
0 @>>> \mathrm H_2(K) @>>> \mathrm H_2(L) @>>> \mathrm H_2(L,K) @>{\mathrm H\partial}>> \mathrm H_1( K) @>>> 0
\end{CD}
\end{equation}
of the pair $(L,K)$ 
shows that the canonical map $\mathrm H_2(K) \to \mathrm H_2(L)$
is an isomorphism if and only if 
the homology boundary 
$\mathrm H_2(L,K) \to \mathrm H_1(K) \cong \mathrm H_1(G)$
is an isomorphism. 

The requirement
that the induced map $C_{\mathrm{Ab}} \to G_{\mathrm{Ab}}$
be an isomorphism of abelian groups
says precisely that the homology boundary 
$\mathrm H_2(L,K) \to \mathrm H_1(K) \cong \mathrm H_1(G)$
is an isomorphism. 
Consequently
the exact homotopy sequence \eqref{central22}
induces the asserted central extension \eqref{central8}.
\end{proof}

\begin{Example}\label{abelian}
Let $A$ be the free abelian group on a set $S$.
Then $\mathrm H_2(A)$ is the exterior square $A \wedge A$
and the free crossed $A$-module $C$ on $S$ fits into a central extension
\begin{equation}
\begin{CD}
0 @>>> A \wedge A @>>> C @>{ \partial}>> A @>>> 1.
\end{CD}
\label{central99}
\end{equation}
This extension is the universal central extension of $A$ that corresponds to
the universal  $A \wedge A$-valued skew-symmetric {\rm 2-\/}cocycle 
on $A$; the skew-symmetricity
determines this {\rm 2-\/}cocycle uniquely.
 In Section {\rm {\ref{univ}}} we will characterize
this kind of central extension by a univeral property.
\end{Example}

\section{Artin's braid group}

The braid group $B_n$ on $n$ strings has $n-1$ generators
$\sigma_1$,\ldots , $\sigma_{n-1}$ subject to the relations
\begin{equation}
\sigma_i\sigma_j=\sigma_j\sigma_i,\ 
\vert i-j\vert \geq 2,
\label{rel1}
\end{equation}
and
\begin{equation}
\sigma_j\sigma_{j+1} \sigma_j= \sigma_{j+1} \sigma_j\sigma_{j+1},
\ 1 \leq j \leq n-2.
\label{rel2}
\end{equation}
Instead of \eqref{rel2}, we take 
the relations
\begin{equation}
\sigma_{j}\sigma_{j+1}\sigma_j \sigma_{j+1}^{-1}\sigma^{-1}_{j}
\sigma_{j+1}^{-1}=e.
\label{rel222}
\end{equation}
Adding the single relation $\sigma_1=e$
plainly kills the entire group, that is, the normal closure of $\sigma_1$
(or any other $\sigma_j$) in $B_n$ is all of $B_n$.

For $1 \leq j \leq n-2$, let
\begin{equation*}
r_{j,j+1}=x_{j}x_{j+1}x_j x_{j+1}^{-1}x^{-1}_{j}
x_{j+1}^{-1}
\end{equation*}
and let
$\widehat R= \{ r_{1,2},r_{2,3},\ldots, r_{n-2,n-1}\}$.
For $ 1 \leq j,k \leq n-1$ with $j < k-1$, 
let
\[
r_{j,k}=[x_j,x_k]=x_jx_k x_j^{-1}x_k^{-1},
\]
let
$\widetilde R$ denote the finite set of these $r_{j,k}$, and
let $R= \widehat R \cup \widetilde R$. Then
\begin{equation}
\langle X, R\rangle  = \langle 
x_1,x_2, \ldots, x_{n-1}; r_{1,2},r_{2,3},\ldots, r_{n-2,n-1},r_{1,3},
\ldots, r_{n-3,n-1}
\rangle
\label{presen}
\end{equation}
is a presentation of $B_n$.
We suppose the notation being adjusted in such a way that,
for $1 \leq j \leq n-1$, the generator $\sigma_j$ 
of  $B_n$ is represented by $x_j$.
Adding the relator $r=x_1$
or, more generally, a relator of the kind $r=x_j$, $1 \leq j \leq n-1$,
we obtain the presentation
$\langle X, R\cup\{r\}\rangle$ of the trivial group.

Let $ K=K(X, R)$ be the geometric 
realization of the presentation
$\langle X, R\rangle$ and let $ L= K\cup e^2$ be the geometric realization
of  the presentation
$\langle X, R\cup\{r\}\rangle$ of the trivial group.
By construction, $\pi_1( K)$ is the group  $ B_n$.
In view of Whitehead's
theorem 
 (Theorem \ref{whit} above),  the homotopy boundary map
$\partial \colon \pi_2( L, K) \to  B_n$
is a free crossed module, 
the $ B_n$-crossed module $\pi_2( L, K)$
being free on the single generator $[r]$.
We will write $\pi_2( L, K)$ as $C_n$.

\begin{Theorem}\label{second}
The kernel of $\partial \colon C_n \to  B_n$
is a central subgroup and
is canonically isomorphic to
$\mathrm H_2(B_n)$. 
Furthermore, for $n=3$,
the group $\mathrm H_2(B_n)$ is zero and for $n \geq 4$,
the group $\mathrm H_2(B_n)$ is cyclic of order {\rm 2},
generated as a subgroup of  $C_n$
by
\begin{equation}
{}^{\sigma_3}[r]
\cdot [r]^{-1}
\in C_n.
\label{I22}
\end{equation}
Consequently, as a crossed module over itself,
$B_n$ is generated by $\sigma_1$, subject to the relation
\begin{equation}
{}^{\sigma_3}\sigma_1=\sigma_1.
\label{rel4}
\end{equation}
\end{Theorem}

Thus the crossed module property corresponds precisely to the relations
\eqref{rel2}, that is, these relations amount to {\em Peiffer identities\/}, 
and the relation \eqref{rel4} correspond exactly to the relations \eqref{rel1}.

N.B. As a crossed module over itself, $B_n$  
has a single generator, which can be taken to be
any of the generators
$\sigma_1$,\ldots , $\sigma_{n-1}$.
In the statement of the theorem, we have written out the relation
\eqref{rel4} with respect to the generator $\sigma_1$.
But we could have written out an equivalent relation
in terms of any of the
 generators
$\sigma_1$,\ldots , $\sigma_{n-1}$.

\begin{Remark}\label{crack} In the language of $B_n$-racks, cf. Example
{\rm \ref{rack}} above, the universal object associated with the
$B_n$-rack $B_n \times \{\sigma_1\}$ (with the obvious structure)
is the free $B_n$-crossed module $C_n$ on the single generator $\sigma_1$.
Thus the universal object associated with the $B_n$-rack 
$B_n \times \{\sigma_1\}$ recovers $B_n$ itself up to the central extension
by $\mathrm H_2(B_n)$. We shall show in Section
{\rm \ref{univ}} below that this central extension is actually universal 
in a suitable sense.
\end{Remark}

\begin{proof} The first statement is a special case of Theorem \ref{basic}.
Indeed, $B_n$ is the normal closure of any of the generators $\sigma_j$,
and 
\[
(C_n)_{\mathrm{Ab}}\cong\pi_2(L,K)_{\mathrm{Ab}}
\cong \mathrm H_2(L,K)\to \mathrm H_1(K)\cong
\mathrm H_1(B_n) \cong \mathbb Z
\]
is plainly an isomorphism.

The group $B_2$ is infinite cyclic and hence has zero second homology group.
The group $B_3$ is well known to be isomorphic to
the trefoil knot group and hence has zero second homology group 
as well.
Henceforth we suppose that $n \geq 4$.

We will now derive the relation \eqref{rel4}.
As usual we denote the 1-skeleton of the cell complex $K$ by $K^1$.
The exact homotopy sequence of the pair $(L,K^1)$ 
takes the form
\begin{equation}
\begin{CD}
0 @>>> \pi_2(L) @>>> \pi_2(L,K^1) 
@>{\partial_L}>> \pi_1(K^1) @>>> 1.
\label{ext1}
\end{CD}
\end{equation}
Since the  $\pi_1(K^1)$-action on $\pi_2(L)$ factors through
$\pi_1(L)$ which is the trivial group,
 the  $\pi_1(K^1)$-action on $\pi_2(L)$
is trivial, and the group extension \eqref{ext1} is necessarily central.
The second homotopy group 
$ \pi_2(L)$ lies in $\pi_2(L,K^1)$ as the group of (classes of)
{\em identities among relations\/} 
in the sense of  \cite{peiffone} and \cite{reideone}.

Let $F\cong \pi_1(K^1)$ denote the free group on $X$. Keeping in mind
that $\pi_2(L,K^1)$ is the free $F$-crossed module on
$R\cup\{r\}$, let
\begin{equation}
\mathcal I_{1,3}=((x_3,r),(e,r^{-1}),(e,r_{1,3}^{-1})).
\label{I1}
\end{equation}
Inspection shows that
$\partial_L([\mathcal I_{1,3}])= e\in F \cong \pi_1(K^1)$ whence
$
[\mathcal I_{1,3}] \in \pi_2(L)\subseteq \pi_2(L,K^1).
$
Indeed, $\mathcal I_{1,3}$ is plainly
an identity among relations.

Consider the canonical projection from
$ \pi_2(L,K^1)$ to $C_n \cong  \pi_2(L,K) $.
Since under this projection 
each generator in $R$ goes to $e \in \pi_2(L,K)$,
the image $[\mathcal I_{1,3}] \in \pi_2(L,K)$
of the member $[\mathcal I_{1,3}]$ of $\pi_2(L)\subseteq \pi_2(L,K^1)$ 
can plainly be written as
\begin{equation}
[\mathcal I_{1,3}]=
{}^{\sigma_3}[r]
\cdot [r]^{-1}
\in \pi_2(L,K)\cong C_n.
\label{I2}
\end{equation}
Here the bracket notation $[\,\cdot \, ]$ is slightly abused.
Since the 2-complex $L$ is simply connected, the Hurewicz map
$\pi_2(L) \to \mathrm H_2(L)$ is an isomorphism.
Furthermore the reasoning in the proof of Theorem \ref{basic}
shows that the canonical map $\mathrm H_2(K) \to \mathrm H_2(L)$
is an isomorphism. Hence, in the exact homotopy sequence
\eqref{ext1}, we can replace $\pi_2(L)$ with  $\mathrm H_2(K)$.
Inspection of the commutative diagram
\begin{equation}
\begin{CD}
0 @>>> \mathrm H_2(K) @>>> \pi_2(L,K^1) 
@>{\partial_L}>> F @>>> 1
\\
@.
@VVV
@VVV
@VVV
@.
\\
0 @>>> \mathrm H_2(B_n) @>>> C_n 
@>{\partial}>> B_n @>>> 1
\label{ext11}
\end{CD}
\end{equation}
reveals that
$[\mathcal I_{1,3}] \in \mathrm H_2(B_n)$.

Likewise, 
for $ 1 \leq j,k \leq n-1$ with $k-j\geq 2$, 
corresponding
 to the relator 
\[
r_{j,k}=[x_j,x_k]=x_jx_k x_j^{-1}x_k^{-1},
\]
we can construct an identity among relations
 $\mathcal I_{j,k}$  
such that $\pi_2(L)\subseteq \pi_2(L,K^1)$
is the free abelian group in 
\[
[\mathcal I_{1,3}],\ldots, [\mathcal I_{n-3,n-1}] \in  \pi_2(L) \subseteq \pi_2(L,K^1).
\]
We shall not
need the explicit form of these identities among relations.
For intelligibility we note that the rank of the free abelian group
$\pi_2(L)\cong \mathrm H_2(L)\cong \mathrm H_2(K)$ equals the cardinality of
the set $\widetilde R$. Indeed, let $\widehat K=K(X,\widehat R)$ be the 
geometric realization of the presentation $\langle X,\widehat R\rangle$.
The canonical map $\mathrm H_2(K) \to \mathrm H_2(K,\widehat K)$
is easily seen to be an isomorphism of abelian groups, the target group
$\mathrm H_2(K,\widehat K)$ being free abelian in the set 
of 2-cells in $K \setminus \widehat K$ or, equivalently, in the set 
$\widetilde R$. 
The association 
\[
\mathcal I_{j,k} \longmapsto r_{j,k}
\quad
(1 \leq j,k \leq n-1,\ j < k-1) 
\]
induces the composite
of the isomorphism $\pi_2(L) \to \mathrm H_2(K)$ with the isomorphism
from $\mathrm H_2(K)$ to $\mathrm H_2(K,\widehat K)$.

For $n \geq 4$, the group
$\mathrm H_2(B_n)$ is well known to be cyclic of order 2,
cf. \cite{cohlamay}, and the left-hand map 
$\mathrm H_2(K) \to \mathrm H_2(B_n)$ in the diagram \eqref{ext11}
is surjective by construction.
For $n=4$, 
the group
$\mathrm H_2(K) \cong \mathbb Z$ is free abelian in the relator $r_{1,3}$,
and
that map takes the form of a surjection
\begin{equation}
\mathbb Z \longrightarrow \mathrm H_2(B_4),
\quad r_{1,3} \longmapsto [\mathcal I_{1,3}].
\label{surj}
\end{equation}
Consequently
$ [\mathcal I_{1,3}] \in  \mathrm H_2(B_4)$ is non-zero and,
the obvious map $\mathrm H_2(B_4) \to \mathrm H_2(B_n)$ being an isomorphism,
we conclude that
$ [\mathcal I_{1,3}] \in  \mathrm H_2(B_n)$ is non-zero for any $n \geq 4$.

Finally, the explicit description \eqref{I2} 
of the element $[\mathcal I_{1,3}] \in C_n$ 
justifies at once the relation \eqref{rel4}
in Theorem \ref{second}.
These observations complete the proof of Theorem \ref{second}.
\end{proof}

\begin{Remark}\label{fcoh}
{\rm The subgroup of $B_n$  ($n \geq 4$) 
generated by $x_1$ and $x_3$ is well known to be
free abelian of rank 
{\rm 2\/}; it is convenient
to think of this group as the direct product $B_2 \times B_2$, viewed
as a subgroup of $B_n$. 
The geometric realization 
$K(x_1,x_3; r_{1,3})$ of the presentation  
$\langle x_1,x_3; r_{1,3}\rangle$ of $B_2\times B_2$ 
is an ordinary  torus which is
plainly aspherical.
Furthermore, 
$\langle x_1,x_3; x_1, r_{1,3}\rangle$ is then a 
presentation of the free cyclic group $C$ generated by $x_3$;
the geometric realization $K(x_1,x_3; x_1,r_{1,3})$ thereof
is a spine for a solid torus having fundamental group $C$,
the identity among relations $\mathcal I_{1,3}$
is defined over the presentation $\langle x_1,x_3; x_1, r_{1,3}\rangle$
of $C$ and, indeed, records the attaching map for a {\rm 3-\/}cell 
to form the solid torus
from the spine. Playing the same game as in the previous proof, but with
the geometric realizations
$K(x_1,x_3; r_{1,3})$ and $K(x_1,x_3; x_1,r_{1,3})$
rather than with $K$ and $L$, respectively, we find
that $\pi_2(K( x_1,x_3; x_1,r_{1,3}))$
is the free $C$-module freely generated by $[\mathcal I_{1,3}]$;
the induced map
from 
$\mathrm H_2(K(x_1,x_3;r_{1,3}))$ to $\mathrm H_2(K(x_1,x_3; x_1,r_{1,3}))$
is an isomorphism, both groups are infinite cyclic,
the group $\mathrm H_2(K(x_1,x_3;r_{1,3}))$ amounts to the second
homology group $\mathrm H_2(B_2\times B_2) \cong \mathbb Z$ of $B_2\times B_2$,
and under the Hurewicz map
from $\pi_2(K( x_1,x_3; x_1,r_{1,3}))$
to $\mathrm H_2(K(x_1,x_3; x_1,r_{1,3}))$ the free $C$-module generator
$[\mathcal I_{1,3}]$ goes to a generator of 
$\mathrm H_2(K(x_1,x_3; x_1,r_{1,3}))$.
Thus the surjection {\rm \eqref{surj}} comes down to the map
\begin{equation}
\mathrm H_2(B_2\times B_2) \longrightarrow \mathrm H_2(B_4) 
\label{inducedh2}
\end{equation}
induced by the injection of $B_2\times B_2$ into $B_4$.
Hence
the map \eqref{inducedh2} yields the generator of
$\mathrm H_2(B_4)$ and, more generally,
of $\mathrm H_2(B_n)$ for $n \geq 4$.
This  observation goes back at least to {\rm \cite{cohlamay}}
where it is pointed out that this kind of construction,
suitably extended, yields the entire integral homology of $B_n$.

The universal covering space of
$K(x_1,x_3; x_1,r_{1,3})$
is an infinite cylinder 
together with infinitely many 2-disks
attached to the cylinder, more precisely, when
the cylinder is realized in ordinary {\rm 3-\/}space,
the disks are attached to the interior of the cylinder.
Thus the identity among relations $\mathcal I_{1,3}$
over the presentation $\langle x_1,x_3; x_1, r_{1,3}\rangle$
of the free cyclic group $C$ is realized geometrically as a 
tesselated sphere arising from the boundary surface of a solid cylinder
having top and bottom a disk labelled $r=x_1$ and wall (lateral surface)
labelled $r_{1,3}$; the cyclinder wall arises from glueing
a rectangle labelled $r_{1,3}$ along the edge labelled $x_3$,
and the sphere arises from glueing the top and bottom disks
along the circles  labelled $x_1$. 
Indeed,  the universal covering space 
of the spine $K(x_1,x_3; x_1, r_{1,3})$
arises from glueing infinitely many copies of that
tesselated sphere along the top and bottom disks.
That tesselated sphere can be seen as a
cycle representing the generator of $\mathrm H_2(B_n)$.}
\end{Remark}

\section{Parallelohedral geometry and 
identities among relations}\label{geometrici}

Using the language of identities among relations,
we will now construct,
under the circumstances of Theorem \ref{second},
a complete system of
$\pi_1(K)$-module generators of $\pi_2(K)$
and thereafter describe the induced map 
$\pi_2(K) \to \pi_2(L)\cong\mathrm H_2(K)$ so that the cokernel
$\mathrm H_2(B_n)\cong \mathbb Z/2$ thereof becomes explicitly visible.
We will in particular exploit the geometry that underlies the hyperplane
arrangements associated with various root systems.
This will provide additional geometric insight into the braid group $B_n$
and will in particular explain {\em why\/} the group $\mathrm H_2(B_n)$ 
is cyclic of order 2 and {\em how\/} the relation \eqref{rel4} arises.
It will also make our reasoning self-contained.
For ease of the reader, we note that the method of {\em pictures\/} developed in
\cite{browhueb}, \cite{collhueb} and \cite{aspheric}
yields a straightforward procedure for reading off the identity
sequences \eqref{J1}, \eqref{J2}, \eqref{J3} and \eqref{J4} below.

We begin with recalling some polyhedral geometry.
A {\em parallelohedron\/} is a convex polyhedron
whose translates by a lattice have disjoint interior and cover
$n$-dimensional real affine space. This notion is implicit in \cite{minkowskione} and was introduced in \cite{voronoione} and
\cite{voronoitwo}.
In dimension 2, the only parallelohedra are the parallelogram and the hexagon.
In dimension 3, 
there are exactly five parallelohedra:  the cube, hexagonal prism, 
truncated octahedron, rhombic dodecahedron, and the 
elongated rhombic dodecahedron. 
As a cellular complex, the 
{\em permutahedron\/} $\mathcal P_n$ of order $n$ is the geometric realization
of the poset of ordered partitions of the set $\{1,2,\ldots,n\}$.
Equivalently, when the symmetric group $S_n$ on $n$ letters
acts on $\mathbb R^n$ via permutation of coordinates,
the permutahedron is the convex hull of the $S_n$-orbit of a generic point.
Thus $S_n$ acts on the permutahedron of order $n$, 
the action being free on the set of vertices. In particular, 
the 1-skeleton of the permutahedron of order $n$ is the Cayley complex of
$S_n$ relative to the standard generators; here the convention is that
any generator that is an involution is represented by a single unoriented edge
rather than by a pair of oppositely oriented edges.
Permutahedra were explored already in \cite{schoutbo}
(not under this name which seems to have been introduced much later).
A permutahedron is a special kind of parallelohedron.
Permutahedra are explored e.~g. on p. 65/66 of \cite{coxemose}.
The permutahedra form an operad \cite{cbergone}, \cite{marsnist} (p.~98).
This fact is behind the cellular models for iterated loop spaces
developed in \cite{milgrone}.
The permutahedron $\mathcal P_3$ is a hexagon and
the permutahedron $\mathcal P_4$ is the truncated
octahedron, that is, 
the Archimedean solid having as faces eight hexagons and six squares,
the corresponding tesselation of
real affine 3-space being the bitruncated cubic honeycomb.
In general, $\mathcal P_n$ is an omnitruncated $(n-1)$-cell.

Given a general hyperplane arrangement,
we will refer to the convex hull $\mathcal A$ 
of a Weyl group orbit of a point in the interior
of a chamber as a {\em Coxeter\/} polytope;
when the Weyl group is finite,
a classifying space for the corresponding Artin group
arises from $\mathcal A$ by the operation of
suitably identifying faces of the polytope $\mathcal A$.
For the braid group, this construction was written down in  
\cite{carlmilg}; the resulting cell complex
is actually dual to the cell complex constructed 
in \cite{foxneuwi} and shown there to be aspherical.  This latter cell  
complex, in turn, arises from a kind of lexicographic cell decomposition
of the configuration space of $n$ particles moving in the plane
yielding a poset that provides a classifying space
for the pure braid group; the operation of taking orbits 
relative to the symmetric group on $n$ letters
then furnishes a classifying space for the ordinary braid group.
For a general Artin group, the classifying space
 arising from a Coxeter polytope
can be found in \cite{salveone} and is also lurking behind the 
approach in \cite{squieron}. 
For the braid group, a suitable comparison
of various operads of the same kind as the permutahedra operad
mentioned above also leads to this classifying space \cite{cbergone}.
The true reason for
the quotient $K$ of the Coxeter polytope
with faces suitably identified
to have the homotopy groups $\pi_j$ zero for $j \geq 2$ 
is that, when $n \geq 2$ is fixed, the 
universal covering space of $K$ is appropriately tiled
by translates of the permutahedron $P_n$;
indeed this universal covering space is a product
of $n-1$ trees (a copy of the real line for $n=2$).

For example, for $j=1,2,3$, the convex hull $\mathcal A_j$ 
of a point in the interior of a chamber
of the root system $A_j$ is an interval for $j=1$,
a hexagon for $j=2$,
and a truncated octahedron for $j=3$. 
Thus $\mathcal A_n$ is exactly the permutahedron $\mathcal P_{n+1}$.
Henceforth we will use the notation  $\mathcal A_n$.
The braid group $B_n$ is the Artin group associated
with the Coxeter graph $A_{n-1}$. In this case, the 
faces are precisely polytopes associated with the full subgraphs
of $A_{n-1}$ with $n-2$ nodes, and
the operation of
suitably identifying faces of the Coxeter polytope
yields a cell complex
whose 2-skeleton is the geometric realization 
$K(X,R)$ of the Artin presentation \eqref{presen}
of $B_n$.

We have already noted that the association 
\[
\mathcal I_{j,k} \longmapsto r_{j,k}
\quad
(1 \leq j,k \leq n-1,\ j < k-1) 
\]
induces the composite
\begin{equation}
\pi_2(L) \longrightarrow 
\mathrm H_2(K,\widehat K)
\label{iso}
\end{equation}
of the isomorphism $\pi_2(L) \to \mathrm H_2(K)$ with the isomorphism
from $\mathrm H_2(K)$ to $\mathrm H_2(K,\widehat K)$.

\begin{Proposition}\label{ident}
As a $B_n$-module, $\pi_2(K)$ has a complete system of generators
\[
[\mathcal J_{j,k,\ell}]\in \pi_2(K)\subseteq \pi_2(K,K^1),
\ 1 \leq j<k<\ell \leq n-1
\]
given by suitable identities $\mathcal J_{j,k,\ell}$ among
the relations of the Artin presentation {\rm \eqref{presen}} of $B_n$.
Furthermore, the composite
of the induced map from $\pi_2(K)$ to $\pi_2(L)$
with  \eqref{iso} 
is given by the associations
\[
[\mathcal J_{j,j+1,j+2}] \longmapsto 2r_{j,j+2},
\
[\mathcal J_{j,j+1,k}] \longmapsto r_{j+1,k} -r_{j,k},
\
[\mathcal J_{j,k-1,k}] \longmapsto r_{j,k} -r_{j,k-1},
\]
in particular, that map is zero on generators of the kind
$[\mathcal J_{j,k,\ell}]$ with $k-j\geq 2$ and $\ell-k\geq 2$.
Consequently $\mathrm H_2(B_n)$ is cyclic of order {\rm 2\/}
and each free generator of $\mathrm H_2(K)$
corresponding to a generator of $\mathrm H_2(K,\widehat K)$ of the kind
$r_{j,k}$ for $k-j\geq 2$ goes to the single non-zero element of
$\mathrm H_2(B_n)$.
\end{Proposition}

We begin with the preparations for
a formal proof of this proposition:
Let $n \geq 4$ and 
consider the Artin presentation 
\eqref{presen} of $B_n$. 
The possible 3-dimensional parallelohedra
that can arise from the full subgraphs of $A_{n-1}$ with three nodes
are the truncated octahedron (Coxeter graph $A_3$),
the hexagonal prism (Coxeter graph $A_1\times A_2$),
and the cube (Coxeter graph $A_1\times A_1\times A_1$).
We will now realize each of these parallelohedra
by an identity among relations.
We will proceed for general $n \geq 4$:

To realize the truncated octahedron we note that,
for each $1 \leq j \leq n-3$, 
the subgroup of $B_n$ generated by $x_j,x_{j+1},x_{j+2}$
is a copy of $B_4$. Accordingly,
consider the  identity 
\begin{equation}
\mathcal J_{j,j+1,j+2} =(q_1,q_2,q_3,q_4,q_5, \ldots, q_{13},q_{14})
\label{J1}
\end{equation}
among relations
arising from reading off an identity sequence from a truncated
octahedron, 
each hexagon being labelled with the relator  
$r_{j,j+1}$ or $r_{j+1,j+2}$ 
and each square being labelled with  the relator $r_{j,j+2}$:
\begin{equation}
\begin{aligned}
q_1&=
(x^{-1}_jx^{-1}_{j+1}x_{j+2}^{-1}x_j^{-1},r_{j,{j+2}})
\\
q_{2}&=
(x_j^{-1}x_{j+1}^{-1}x_{j+2}^{-1}x_j^{-1}x_{j+2},r_{j,j+1})
\\
q_{3}&=
(x_j^{-1}x_{j+1}^{-1}x_{j+2}^{-1}x_j^{-1}x_{j+1}^{-1},r_{j+1,j+2})
\\
q_{4}&=
(x_j^{-1}x_{j+1}^{-1}x_{j+2}^{-1}x_j^{-1}x_{j+1}^{-1}x_{j+2}x_{j+1},r_{j,{j+2}}^{-1})
\\
q_{5}&=
(x_j^{-1}x_{j+1}^{-1}x_{j+2}^{-1}x_{j+1}^{-1},r_{j+1,j+2}^{-1})
\\
q_{6}&=
(x_j^{-1}x_{j+2}^{-1}x_{j+1}^{-1}x_{j+2}^{-1}x_j^{-1}x_{j+1}^{-1},r_{j,j+1}^{-1})
\\
q_{7}&=
(x_j^{-1}x_{j+2}^{-1}x_{j+1}^{-1}x_{j+2}^{-1}x_j^{-1}x_{j+1}^{-1},r_{j,{j+2}})
\\
q_{8}&=
(x_j^{-1}x_{j+2}^{-1}x_{j+1}^{-1}x_{j+2}^{-1}x_j^{-1}x_{j+1}^{-1}x_{j+2},r_{j,j+1})
\\
q_{9}&=
(x_j^{-1}x_{j+2}^{-1}x_{j+1}^{-1}x_{j+2}^{-1}x_j^{-1}x_{j+1}^{-1}x_{j+2}x_{j+1}x_j,r_{j+1,j+2})
\\
q_{10}&=
(x_j^{-1}x_{j+2}^{-1}x_{j+1}^{-1}x_{j+2}^{-1}x_j^{-1},r_{j,{j+2}}^{-1})
\\
q_{11}&=
(x_j^{-1}x_{j+2}^{-1}x_{j+1}^{-1}x_j^{-1}x_{j+2}^{-1}x_{j+1}^{-1},r_{j+1,j+2}^{-1})
\\
q_{12}&=
(x_j^{-1}x_{j+2}^{-1}x_{j+1}^{-1}x_j^{-1} x_{j+1} x_{j+2}^{-1},r_{j,{j+2}})
\\
q_{13}&=
(x_j^{-1}x_{j+2}^{-1}x_{j+1}^{-1}x_j^{-1} ,r_{j,j+1}^{-1})
\\
q_{14}&=
(x_{j+2}^{-1}x_j^{-1},r_{j,{j+2}})
\end{aligned}
\end{equation}
For $B_4$, the construction stops at this stage.
In particular, attaching a 3-cell to 
the geometric realization
$K(x_1,x_2,x_3; r_{1,2},r_{2,3},r_{1,3})$
via the attaching map given by $\mathcal I_{1,2,3}$
yields a classifying space of $B_4$.
This amounts precisely the operation of
identifying faces of the truncated octahedron mentioned earlier.
Thus the kernel of 
$\pi_2(K,K^1)_{\mathrm {Ab}}\to \mathbb Z[B_4]\langle x_1,x_2,x_3\rangle $
is the free $B_4$-module freely generated by $[\mathcal I_{1,2,3}]$ and
the resulting $B_4$-chain complex is a free resolution of the integers
over $B_4$; cf. e.~g. \cite{squieron}. Consequently
\begin{equation}
\begin{CD}
0 @>>> \pi_2(K) @>>> \pi_2(K,K^1) 
@>{\partial}>> \pi_1(K^1) @>>> B_4 @>>> 1
\end{CD}
\end{equation}
is then a free crossed resolution of $B_4$; cf. \cite{crossed}
for the notion of free crossed resolution.

When $n \geq 5$, additional identities among relations
arise, according to the possible full subgraphs $A_2\times A_1$ of 
$A_{n-1}$:
For each pair $(j,k)$ with
$1 \leq j<k \leq n-1$ and  $j < k-2$,
the subgroup of $B_n$ generated by $x_j,x_{j+1},x_k$
is isomorphic to the direct product 
$B_3 \times B_2$
of a copy of  $B_3$ with a copy of $B_2$.
Accordingly,
consider
the identity 
\begin{equation}
\mathcal J_{j,j+1,k} =(q_1,q_2,q_3,q_4,q_5, q_6,q_7,q_8)
\label{J2}
\end{equation}
among relations
arising from reading off an identity sequence from a prism 
having base and top a hexagon labelled $r_{j,j+1}$ and
having, furthermore, six lateral squares,
each square being labelled with  the relators $r_{j,k}$ or $r_{j+1,k}$:
\begin{equation}
\begin{aligned}
q_1&=(x_{j+1}^{-1}x_j^{-1}x_{j+1}^{-1},r_{j,j+1}^{-1})
\\
q_2&=(x_{j+1}^{-1}x_j^{-1}x_{j+1}^{-1}x_jx_k^{-1},r_{j+1,k}^{-1})
\\
q_3&=(x_{j+1}^{-1}x_j^{-1}x_{j+1}^{-1}x_k^{-1},r_{j,k}^{-1})
\\
q_4&=(x_{j+1}^{-1}x_j^{-1}x_{j+1}^{-1}x_k^{-1},r_{j+1,k})
\\
q_5&=(x_{j+1}^{-1}x_j^{-1}x_k^{-1}x_{j+1}^{-1}x_jx_{j+1},r_{j,k}^{-1})
\\
q_6&=(x_{j+1}^{-1}x_j^{-1}x_k^{-1}x_{j+1}^{-1},r_{j,j+1})
\\
q_7&=(x_{j+1}^{-1}x_j^{-1}x_k^{-1},r_{j,k})
\\
q_8&=(x_{j+1}^{-1}x_k^{-1},r_{j+1,k}) .
\end{aligned}
\end{equation}
With the notation $F_{3,2}$ for the free
group on $x_j,x_{j+1}, x_k$ and $C_{3,2}$ for the free crossed
$F_{3,2}$-module on the relators $r_{j,j+1},r_{j,k},r_{j+1,k}$, the sequence
\begin{equation}
\begin{CD}
0 @>>> \mathbb Z[B_3\times B_2]\langle \mathcal J_{j,j+1,k}\rangle 
@>>> C_{3,2} 
@>{\partial}>> F_{3,2} @>>> B_3 \times B_2 @>>> 1
\end{CD}
\end{equation}
is a free crossed resolution of $B_3\times B_2$.

Likewise, for each pair $(j,k)$ with
 $1 \leq j<k \leq n-1$ and  $j < k-2$,
the subgroup of $B_n$ generated by $x_j,x_{k-1},x_k$
is isomorphic to the direct product
$B_2 \times B_3$
of a copy of $B_2$ with a copy of $B_3$.
Accordingly,
consider
the identity 
\begin{equation}
\mathcal J_{j,k-1,k} =(q_1,q_2,q_3,q_4,q_5, q_6,q_7,q_8)
\label{J3}
\end{equation}
among relations
arising from reading off an identity sequence from a prism 
having base and top a hexagon labelled $r_{k-1,k}$ and
having, furthermore, six lateral squares,
each square being labelled with  the relators $r_{j,k-1}$ or $r_{j,k}$:
\begin{equation}
\begin{aligned}
q_1&=(x_k^{-1}x_{k-1}^{-1}x_k^{-1},r_{k-1,k}^{-1})
\\
q_2&=(x_k^{-1}x_{k-1}^{-1}x_k^{-1}x_{k-1}x_j^{-1},r_{j,k})
\\
q_3&=(x_k^{-1}x_{k-1}^{-1}x_k^{-1}x_j^{-1},r_{j,k-1})
\\
q_4&=(x_k^{-1}x_{k-1}^{-1}x_k^{-1}x_j^{-1},r_{j,k}^{-1})
\\
q_5&=(x_k^{-1}x_{k-1}^{-1}x_j^{-1}x_k^{-1}x_{k-1}x_k,r_{j,k-1})
\\
q_6&=(x_k^{-1}x_{k-1}^{-1}x_j^{-1}x_k^{-1},r_{k-1,k})
\\
q_7&=(x_k^{-1}x_{k-1}^{-1}x_j^{-1},r_{j,k-1}^{-1})
\\
q_8&=(x_k^{-1}x_j^{-1},r_{j,k}^{-1})
\end{aligned}
\end{equation}

When $n \geq 6$, more identities among relations arise,
according to the possible full subgraphs $A_1\times A_1\times A_1$ 
of $A_{n-1}$:
For each triple $(j,k,\ell)$ with
$1 \leq j<k< \ell \leq n-1$, $k-j\geq 2$, and $\ell -k \geq 2$,
the subgroup of $B_n$ generated by $x_j,x_k,x_{\ell}$
is free abelian of rank three.
The familiar commutator identity
\[
[x,y]\cdot {}^y[x,z]\cdot [y,z]\cdot {}^z[y,x]\cdot [z,x]\cdot {}^x[z,y]=e
\]
valid in any group entails the following
identity 
\begin{equation}
\mathcal J_{j,k,\ell} =(q_1,q_2,q_3,q_4,q_5, q_6)
=(
(e,r_{i,j}),
(x_j,r_{i,k})
(e,r_{j,k}),
(x_k,r_{i,j}^{-1})
(e,r_{i,k}^{-1}),
(x_i,r_{j,k}^{-1}))
\label{J4}
\end{equation}
among relations.
This identity among relations also arises
from reading off an identity sequence from a cube 
whose pairs of opposite lateral squares 
are labelled with the relators $r_{j,k}$, $r_{j,\ell}$, or $r_{k,\ell}$.
With the notation $F_{2,2,2}$ for the free
group on $x_j,x_k, x_\ell$ and $C_{2,2,2}$ for the free crossed
$F_{2,2,2}$-module on the relators $r_{j,k},r_{j,\ell},r_{k,\ell}$,
the resulting sequence
\begin{equation}
\begin{CD}
0 @>>> \mathbb Z[B_2^{\times 3}]\langle \mathcal J_{j,k,\ell}\rangle 
@>>> C_{2,2,2} 
@>{\partial}>> F_{2,2,2} @>>> B_2^{\times 3} @>>> 1
\end{CD}
\end{equation}
is then a free crossed resolution of $B_2\times B_2\times B_2$.

We can now prove Proposition \ref{ident}.
When $n \geq 4$, as a $B_n$-module, 
$\pi_2(K)$ is generated by 
the elements 
\[
[\mathcal J_{j,k,\ell}]\in \pi_2(K)\subseteq \pi_2(K,K^1),
\ 1 \leq j<k<\ell \leq n-1.
\]
Indeed, the induced map $\pi_2(K) \to \pi_2(K,K^1)_{\mathrm {Ab}}$
is still injective  and $\pi_2(K,K^1)_{\mathrm {Ab}}$ amounts to the free
$B_n$-module freely generated by the 2-cells of $K$ or, equivalently,
by the set $R$ of relators. 
Let $ \mathbb Z[B_n]$ denote the integral group ring of
$B_n$ and let
$ \mathbb Z[B_n]\langle X\rangle $
be  the free  $B_n$-module 
freely generated by $X$.
The  map
$\pi_2(K,K^1)_{\mathrm {Ab}}\to \mathbb Z[B_n]\langle X\rangle $
induced by $\partial \colon \pi_2(K,K^1) \to \pi_1(K^1)$
has kernel the isomorphic image of $\pi_2(K)$
and, together with the obvious map
$\mathbb Z[B_n]\langle X\rangle \to \mathbb Z[B_n]$,
yields the beginning of the  familiar free resolution
of the integers in the category of $ \mathbb Z[B_n]$-modules,
well explored at various places in the literature
\cite{deligone}, \cite{fadeneuw}, \cite{salveone}, \cite{squieron}. 
These remarks establish Proposition \ref{ident}.
The 14-gon and the octagons discussed above are related with the geometry of
posets of maximal chains of the symmetric group
on $n$ letters, cf. e.~g. \cite{lawreone}.

\section{Universal central extensions}\label{univ}

The Schur multiplicator of a group (later identified as the second homology group)
was originally discovered by Schur in the search for central extensions 
to realize projective representations by linear ones.
In this spirit, we will now push further the situation of Example \ref{ex2}.
The starting point is the (familiar)
observation that, given the group $G$
with $\mathrm H_1(G)$ free abelian or zero,
for any abelian group $\AA$,
 the universal 
coefficient map
\begin{equation}
\begin{CD}
\mathrm H^2(G,\AA) @>>>  \mathrm{Hom}(\mathrm H_2(G),\AA)  
\end{CD}
\label{univco}
\end{equation}
is an isomorphism whence
the central extensions of $G$ by the abelian group $\AA$
are then classified by $\mathrm{Hom}(\mathrm H_2(G),\AA)$.

We will refer to a central extension
\begin{equation}
\begin{CD}
0 @>>> \mathrm H_2(G) @>>> U @>>> G @>>> 1
\label{un}
\end{CD}
\end{equation}
as being {\em universal\/} provided it has the following universal property:
Given the central extension 
\begin{equation*}
\begin{CD}
\mathbf e \colon 
0
@>>>
\AA
@>>>
E
@>{\pi}>>
G
@>>>
1,
\end{CD}
\end{equation*}
the identity of $G$ lifts to a commutative diagram
\begin{equation}
\begin{CD}
0 @>>> \mathrm H_2(G) @>>> U @>>> G @>>> 1
\\
@.
@VV{\vartheta}V
@VV{\Theta}V
@VV{\mathrm{Id}}V
@.
\\
0
@>>>
\AA
@>>>
E
@>{\pi}>>
G
@>>>
1 
\label{CD88}
\end{CD}
\end{equation}
of central extensions in such a way that
the assignment to $\mathbf e$ of $\vartheta$ induces the universal  
coefficient map {\rm \eqref{univco}}.

When $\mathrm H_1(G)$ is free abelian or zero (that is, $G$ is perfect), 
it is straightforward to
construct a universal central extension of $G$ and
standard abstract nonsense reveals that a universal extension is 
then unique
up to isomorphism.
Here the notion of universal central extension is to be 
interpreted with a grain of salt, though,
since this notion is lingua franca only for the case where
$G$ is perfect \cite{milnokbo} and can then be characterized in other ways.

\begin{Theorem} Under the circumstances of Theorem {\rm \ref{basic}},
the central extension
{\rm \eqref{central8}}
is the universal central extension of $G$.
\end{Theorem}

\begin{proof} This is an immediate consequence of the freeness
of the $G$-crossed module $C$. We leave the details to the reader.
\end{proof}

Thus, in particular, the free $B_n$-crossed module  $C_n$
in Theorem \ref{second} above yields the universal central extension
of $B_n$. Likewise the central extension
\eqref{central99} is the universal central extension
of the free abelian group $A$ on the set $S$.
Furthermore, the map {\rm \eqref{surj}}
is the universal map induced from the injection of $A$ in $\mathrm B_4$
(or more generally $\mathrm B_n$)
via the universal property of
the universal central extension.

\end{document}